\documentclass{amsart}
\usepackage{amssymb, array, verbatim, amscd}
\usepackage{amsmath, amsfonts,amsthm}

\begin{document}
\title[ Finite Dynamical Systems, Linear Automata, and Finite Fields]{ Finite Dynamical Systems, Linear Automata, and Finite Fields}
\author
{O. MORENO}
\address{ Department of Mathematics and Computer Science\\
 University of Puerto Rico at Rio Piedras\\
 Rio Piedras, PR 00931, moreno@uprr.pr}
 \author{D. BOLLMAN}
 \address{Department of Mathematics\\
 University of Puerto Rico at Mayaguez\\
 Mayaguez, PR 00681-9018\\
 bollman@cs.uprm.edu }
\author{M.A. AVI\~{N}O-DIAZ}
 \address{Department of Mathematics\\
 University of Puerto Rico at Cayey\\
 Cayey,PR 00777\\
 mavino@cayey.upr.edu}

\maketitle
\begin{abstract}
 We establish a connection between finite fields and finite
dynamical systems. We show how this connection can be used to shed light
on some problems in finite dynamical systems and in particular, in
linear systems.
\end{abstract}

\medskip
\medskip
\section{ Introduction}
There is a natural correspondence between the set $GF(p^r)$ and
the set ${{\bf Z}_p}^r$ of  $r$-tuples over ${\bf Z}_p,$ p prime.
Furthermore, $GF(p^r)$ is a vector space over $GF(p)={\bf Z}_p$
and linear transformations over ${{\bf Z}_p}^r$ correspond to
linearized polynomials over $GF(p^r).$ In this ongoing work, we
use these facts to study some problems in finite dynamical systems
and linear automata.
\par
In Section 2, we study finite dynamical systems and how our approach
can be used in  the classification problem. In Section 3, we study
linearized polynomials and how they can be applied to a problem in
linear finite state machines, which in turn arises from a problem in
crystallographic FFTs.
\section{Finite dynamical systems}

Finite dynamical systems are important in applications to
computational molecular biology. They are useful in microarrays of
genes in order to find the best model that fits a given data, the
so called ``reverse engineering problem.'' (See
http://industry.ebi.ac.uk/~brazama/Genenets.) Laubenbacher and
Pareigis [2] define a finite dynamical system  as a function
$f:k^n\rightarrow k^n,$  constructed by the following data:

1. $k=\{0,1\}$

2. a finite graph $F$ on $n$ vertices.

3. a family of ``local" update functions $f_a:k^n\rightarrow k^n$,
one for each vertex $a\in F$, which changes only the coordinate
corresponding to $a$, and computes the binary state of vertex $a$.
These functions are local in the sense that they only depend on
those variables which are connected to $a\in F$.

4. an ``update schedule" $\pi$, which specifies an order on the
vertices of $F$, represented by a permutation $\pi\in S_n$.

 The function $f$ is then constructed by composing the local
 functions according to the update schedule $\pi$, that is
 $$f=f_{\pi (n)} \circ \cdots \circ f_{\pi (1)}: k^n\rightarrow
 k^n.$$

 In [3], Laubenbacher and Pareigis called the above function
  a permutation sequential dynamical system,  and extended
  the definition in
different directions. In particular, they take  $k$ as  an
arbitrary set. For our purposes it is convenient to consider
$k={\bf Z}_m=\{0,1,\cdots,m-1\}$,
\medskip
\par\noindent
Definition 2.1. A finite dynamical system (FDS) is a pair $(V,f)$
where $V$ is the {\it set} of vectors over a finite field and
$f:V\rightarrow V.$
\medskip
\par\noindent
Definition 2.2. The state diagram of a FDS $(V,f)$ is the digraph
whose vertices are members of $V$ and whose edges are the set of
all $(x,f(x)),$ where $x\in V.$
\medskip
\par\noindent
Remark: Note how for an FDS $({{\bf Z}_p}^r,f),$ the function can
be viewed naturally as an $r$-tuple $(f_1,f_2,\cdots,f_r)$ of
functions $f_i:{\bf Z}_p\rightarrow {\bf Z}_p,$ $i\le i \le r.$ It
is also important to note that, using Langrange interpolation, any
function from a finite field to itself can be realized as a
polynomial [4]. Hence each of the $f_i$ can be regarded as a
polynomial over ${\bf Z}_p.$ A similar remark applies to a FDS of
the form $(GF(p^n),f).$
\medskip
\par\noindent
Definition 2.3. Two FDSs are isomorphic if their state diagrams
are isomorphic.
\medskip
\par
Two FDSs are isomorphic if and only if their state diagram  are
isomorphic as digraphs. A \textit{limit cycle} is simply a directed
cycle in the state diagram ${\mathcal S}_f.$ A loop is a limit cycle
consisting of a single vertex and in the case that it occurs, it is
a fixed point of the FDS $f.$ We denote by ${\mathcal L}_f$ the
subdigraph of ${\mathcal S}_f $ induced by all the arcs of the limit
cycles.
\medskip
\par
Definition 2.4. Let $f:V\rightarrow V$ be a FDS with state diagram
${\mathcal S}_f$ and with subdigraph ${\mathcal L}_f$ of limit
cycles. Then $x\in k$ is a vertex in ${\mathcal L}_f$ if and only if
there exists a positive integer $m$ such that
 $f^{m}(x)=x$. The minimum $m$ such that $f^{m}(x)=x$ for all $
x\in {\mathcal L}_f$ is  called  the \textit{order} of the system
$f$, denoted by Order$(f)$. (See [2])
\medskip
\par\noindent
Directed paths in ${\mathcal S}_f$ correspond to iterations of $f$
on the element at the beginning of the path. Since the set ${{\bf
Z}_p}^r$ is finite, any directed path must eventually enter a limit
cycle. Thus each connected component of ${\mathcal S}_f$ consist of
one limit cycle, together with \textit{transients}, that is,
directed paths having no repeated vertices and ending in a vertex
that is part of limit cycle.
\medskip\medskip
\medskip
\par
One of the main problems in FDSs is their classification. Loosely speaking,
this is the problem of determining of two arbitrarily given FDSs whether or
not they are isomorphic. One of our goals
in this work is to facilitate the solution of the classification problem
though the association given in
\medskip
\par\noindent
Theorem 2.5. For any fixed basis $\alpha_1,\cdots,\alpha_r$ of
$GF(p^r)$ there is a natural one-one correspondence between the
FDSs over $GF(p^r)$ and those over ${{\bf Z}_p}^r.$
\medskip
\par\noindent
{\it Proof.} There is a natural correspondence between the sets
${{\bf Z}_p}^r$ and $GF(p^r)$, namely,
$(x_1,x_2,\cdots,x_r)\leftrightarrow x_1\alpha_1+\cdots
x_r\alpha_r.$ Now given $f:{{\bf Z}_p}^r \rightarrow {{\bf
Z}_p}^r,$ define $L:GF(p^r)\rightarrow GF(p^r)$ such that for each
$(x_1,x_2,\cdots,x_r)$ in ${{\bf Z}_p}^r,$ $L(x)=f_1(x_1)\alpha_1
+ \cdots f_r(x_r)\alpha_r,$ where $x=x_1\alpha_1 + \cdots +
x_r\alpha_r$ and $f=(f_1,\cdots,f_r).$ Conversely,  given
$L:GF(p^r)\rightarrow GF(p^r),$ if $L(x_1\alpha_1+\cdots
x_r\alpha_r)=y_1\alpha_1+\cdots y_r\alpha_r,$ there corresponds a
function $f=(f_1,f_2,\cdots,f_r)$ such that
$f(x_1,\cdots,x_r)=(y_1,\cdots,y_r),$ where $f_i(x_i)=y_i$ for
each $i=1,2,\cdots r.$ Since this correspondence is onto and the
two sets are finite with the same number of elements, it is also
one-one.
\medskip
\par\noindent
Corollary 2.6. If $S_1=({{\bf Z}_p}^r,f)$ and $S_2=({{\bf
Z}_p}^r,f')$ are FDSs and $f$ corresponds to $L(x)$ with respect
to the basis $\alpha_1,\cdots,\alpha_r$ and $f'$ corresponds also
to $L(x),$ but with respect to another basis, then $S_1$ and $S_2$
are isomorphic.
\medskip
\par
This latter corollary says that our approach is  quite useful for the
classification problem. On the other hand,
\medskip
\par\noindent
Theorem 2.7. For any fixed basis $\alpha_1,\cdots,\alpha_r$ of
$GF(p^r),$ there is a natural correspondence between the FDSs over
$(GF(p^r))^n$ and those over $({\bf Z}_p)^{rn}.$
\medskip
\par
In other words, it is redundant to study both types of these FDSs,
but each is important given the classification problem.
\section{Linear FDSs and linearized polynomials.}
 A linear finite dynamical
system or a linear (autonomous) finite state machine is a FDS
$({{\bf Z}_p}^r,f)$ in which $f$ is a linear transformation on
${{\bf Z}_p}^r$ regarded as a vector space over ${\bf Z}_p.$ We
shall see in this case that there is a useful correspondence between
linear FDSs and linearized polynomials.
\par
The correspondence $x\rightarrow x^{p^i},$ $i=0,1,\cdots,p^{r-1},$
gives the Galois automorphisms of $GF(p^r).$ A linearized
polynomial $L(x)$ is a polynomial generated by these
automorphisms. In other words, $L(x)=\sum_{i=0}^{r-1}A_ix^{p^i},$
where $A_i \in GF(p^r.)$ We note that if $y,z \in GF(p^r)$ and
$\lambda \in GF(p),$ then $L(x+y)=L(x)+L(y)$ and $L\lambda x) =
\lambda L(x).$ Thus, $L(x)$ is a linear function on $GF(p^r)$
regarded as a vector space over $GF(p).$ Furthermore the
correspondence between $GF(p^r)$ and ${{\bf Z}_p}^r$ given in
Theorem 1 is an isomorphism as a vector space over ${\bf Z}_p$ .
Since there are $(p^r)^r$ linearized polynomials, this coincides
with all the linear functions on ${{\bf Z}_p}^r.$ It is easy to
see that if $f:{{\bf Z}_p}^r \rightarrow {{\bf Z}_p}^r$ is a FDS
associated to the linearized polynomial
$L(x)=\sum_{i=0}^{r-1}A_ix^{p^i},$ then ker$f$ is the set of all
roots of $L(x)$. So,  $f$ is invertible if and only if the only
root of $L$ in $GF(p^r)$ is $0$.
\par
Given a linear autonomous machines $S=({{\bf Z}_p}^r,F),$ if $f$
is a nonsingular linear transformation, i.e., an invertible matrix
over ${\bf Z}_p,$ then the state space ${{\bf Z}_p}^r$ decomposes
into disjoint ``orbits''  or ``cycles.'' Based on the above
observation, the same also holds for machines $(GF(p^r),L)$ and it
is interesting to note how properties of linearized polynomials
determine this orbit structure, much in the same way as the
properties of the elementary divisors of $f$ determine the orbit
structure of $({\bf Z}_p^r,f)$ in the classical theory. However,
our motivation for studying linearized polynomials stems from a
more general problem which arises in crystallographic FFTs [5].
Let us briefly describe this problem.
\par
Crystallographic data can introduce structured symmetries in the
inputs of a multidimensional discrete Fourier transform,   which
in term introduce symmetries into the outputs. In order to avoid
redundant calculations, it is of interest to exploit these
symmetries.  Assuming that symmetries are given by an $n\times n$
matrix $S$ over ${\bf Z}_p,$ for  prime $p$ edge length, we can
reduce the complexity of the FFT by determining a matrix $M$ with
$MS=SM$ and $M^tS=SM^t$ (where $M^t$ denotes the transpose of $M$)
that minimizes the number of ``$MS$-orbits.'' A vector $x\in
Z_p^n$ belongs to an {\it MS-orbit} of length $k$ if and only if
$M^kx=S^ix$ for some $i.$ The cases $n=2$ and $n=3$ are of
particular interest.
\par
For $n=2,$ for example, i.e., $F:{\bf Z}_p\times {\bf Z}_p
\rightarrow {\bf Z}_p\times {\bf Z}_p,$ this corresponds to the
study of linearized polynomials $F(x)=Ax^p+Bx,$ where $A,B\in
GF(p^2),$ where $F$ is an invertible map. The first question is,
when is $F(x)$ invertible in $GF(p^2)?$
\medskip
\par\noindent
Lemma 3.1. $F(x)$ is invertible if $A^{p+1}\ne B^{p+1}.$
\medskip
\par\noindent
{\it Proof.} $F(x)$ is invertible if and only if it is one-one,
i.e., if $ker\;F=0.$
In other words, if $F(x)=0$ has only $x=0$ as a solution over $GF(p^2).$
But $x\ne 0$ and $Ax^p+Bx=0$ imply that
$X^{p-1}=-{{B}\over{A}}$ and so raising both sides to the power $p+1,$
we obtain $x^{p^2-1}={{B^{p+1}}\over{A^{p+1}}}$ and $A^{p+1}=B^{p+1}.$
\medskip
\par
In the remainder of this section we consider the class ${\mathcal
L}_p$ of linearized polynomials $L(x)=\sum_{i=0}^{r-1}A_ix^{p^i}$
where $A_i\in GF(p).$ These types of polynomials have important
properties which we outline below.
\medskip
\par\noindent
Property I. If $L(x),L'(x)\in {\mathcal L}_p$ then
$L(L'(x))=L'(L(x)).$ (Note that this means that the corresponding
matrices commute).
\medskip
\par\noindent
Property II. Given $L(x)\in {\mathcal L}_p,$ the class of $L'(x)$
satisfying Property I is precisely ${\mathcal L}_p.$
\medskip
\par\noindent
Property III. Using a normal basis, the matrix corresponding to
$L(x)\in {\mathcal L}_p$ is symmetric. Consequently, the transpose
matrix also commutes.
\medskip
\par\noindent
Definition 3.2. If $L(x)=\sum_{i=0}^{r-1}A_ix^{p^i}$ is a
linearized polynomial, then its {\it associate} is
$l(x)=\sum_{i=0}^{r-1}A_iX^i.$
\medskip
\par\noindent
Definition 3.3. $L^i(x)=L(x)$ and $L^{n+1}(x)=L(L^n(x)).$ That is,
$L^n(x)$ is the $n$-fold composition of $L$ with itself.
\medskip
\par\noindent
Property IV. $L^n(x)=x$ modulo $x^{p^r}=x$ if and only if
$(l(x))^n=1$ modulo $x^r=1.$
\section{ Systems over ${\bf Z}_{p^n}$}

In this section we study FDS over ${\bf Z}_{p^n}$, that is systems
$(V,f)$ where  $f:{{\bf Z}_{p^n}}^r\rightarrow {{\bf Z}_{p^n}}^r$
and $V$ is the ${\bf Z}_{p^n}$-module ${{\bf Z}_{p^n}}^r$. We use
Definitions 2.2, 2.3, and 2.4,
 but now we are working  in the
${\bf Z}_{p^n}$-module $V$. Therefore a linear FDS is an
endomorphism of the ${\bf Z}_{p^n}$-module $V$.

Let $g:{{\bf Z}_p}^n\rightarrow {\bf Z}_{p^n}$ be  a bijection.
Then the product function $g^r:{{\bf Z}_p}^{nr}\rightarrow {\bf
Z}_{p^n}^r$ given by $g (a_1,a_2,\cdots
,a_r)=(g(a_1),g(a_2),\cdots ,g(a_r))$ is a bijection too.
\medskip
\par\noindent
Proposition 4.1.
  Let $f:{{\bf Z}_{p^n}}^r\rightarrow {{\bf Z}_{p^n}}^r$ be a FDS. Let
  $\overline f:{{\bf Z}_p}^{nr}\rightarrow {{\bf Z}_p}^{nr}$
  be the FDS  such that
 $g^r\circ \overline f= f \circ g^r$. Then
  $\overline f$ and  $f$ have  state diagrams isomorphic.
\medskip
\par\noindent
{\it Proof.}
 Since $g^r$ is a bijection, there exists
  $(g^r)^{-1}=g^{-r}$. Then the system
   $f_1=g^{-r}\circ f \circ
  g^r$ has the same state diagram to $f$. In fact, set $x_1=g^{-r}(x)$ and $y_1=g^{-r}(y)=(g^{-r}\circ f \circ
  g^r)(x_1)$. Now,  suppose $(x,y=f(x))$ is an edge in the state
  diagram of $f$. Then $(x_1, y_1)$ is an edge in the state diagram of $f_1$.
   On the other hand, $g^r\circ  f_1= f \circ
  g^r$. So, $\overline f=f_1$ and our claim holds.
\medskip
\par\noindent
\medskip
Definition 4.2. With the notation above, if $f$ is a linear FDS
over ${\bf Z}_{p^n}$ then the system $f$ will be called the linear
system associated to $\overline f$ by the bijection $g$.
\medskip
\par\noindent
\medskip
We use the the linear system associated to a non-linear  FDS
$\overline f$ to describe its state diagram and its order.

\medskip\par\noindent
Example

  Let $f:({\bf Z}_{2^3})^2\rightarrow ({\bf Z}_{2^3})^2$ be a linear system, given by
$f(a,b)=(5b,a+2b)$. For any bijection $g:({\bf Z}_2)^3\rightarrow
{\bf Z}_{2^3}$ we have an induced systems $\overline f$ with the
same state diagram of the system $f$.  Now using the method given
in [1], we find the Order of $f$ and the Order of $\overline f$
for any bijection $g$. The matrix
$$A=\left(\begin{array}{cc}
  0 & 5 \\
  1 & 2
\end{array}\right)\equiv \left(\begin{array}{cc}
  0& 1 \\
  1 & 0
\end{array}\right) \ (\hbox{ mod $2$ } )$$ has minimal polynomial
 \[m(x)=(x-1)^2\] and the Order of $A$ modulo $2$ is
$e=2$. Since \[A^2=\left(\begin{array}{cc}
  5 & 2 \\
  2 & 1
\end{array}\right) (\hbox{ mod $2^3$ } ),\] the largest  positive integer
$\beta$ such that $A^2\equiv I \ (\hbox{ mod } 2^\beta)$ is $\beta
=1$. Then $A^2$ has Order  $4$ and $A$ has Order $8$ modulo $2^3$.

\section{ Future Work}
We will exploit the ideas presented here to seek a polynomial
solution to the reverse engineering problem. We  also make use of
our theory to develop an efficient algorithm to determine, given a
matrix $S$ of symmetries, a matrix $M$ that minimizes the number of
$MS$-orbits in the precomputation phase of crystallographic FFTs.
\smallskip
\section*{Acknowledgments}
The work of the first two authors
was supported in part by the NSF grants CISE-MI and NIH/NIGMS S06GM08103.

\end{document}